\documentclass[a4paper,12pt]{article}
\usepackage{amscd}
\usepackage{amsmath,amsfonts,amssymb,amscd}
\usepackage{indentfirst,graphicx,epsfig}
\usepackage{graphicx,psfrag}
\input{epsf}

\newtheorem{thm}{Theorem}

\newtheorem{cor}{Corollary}

\newtheorem{lem}{Lemma}

\def\qed{\hfill \nopagebreak\rule{5pt}{8pt}}

\title{\bf\Large A Tur\'an-type problem on degree
sequence\footnote{Supported by NSFC and the ``973" project. } }
\author{
\small  Xueliang Li, Yongtang Shi\footnote{Corresponding author. }\\
\small Center for Combinatorics and LPMC-TJKLC \\
\small Nankai University, Tianjin 300071, China \\
\small lxl@nankai.edu.cn,  shi@nankai.edu.cn
\date{}}
\begin{document}
\maketitle
\begin{abstract}
Given $p\geq 0$ and a graph $G$ whose degree sequence is
$d_1,d_2,\ldots,d_n$, let $e_p(G)=\sum_{i=1}^n d_i^p$. Caro and
Yuster introduced a Tur\'an-type problem for $e_p(G)$: given $p\geq
0$, how large can $e_p(G)$ be if $G$ has no subgraph of a particular
type. Denote by $ex_p(n,H)$ the maximum value of $e_p(G)$ taken over
all graphs with $n$ vertices that do not contain $H$ as a subgraph.
Clearly, $ex_1(n,H)=2ex(n,H)$, where $ex(n,H)$ denotes the classical
Tur\'an number, i.e., the maximum number of edges among all $H$-free
graphs with $n$ vertices. Pikhurko and Taraz generalize this
Tur\'an-type problem: let $f$ be a non-negative increasing real
function and $e_f(G)=\sum_{i=1}^n f(d_i)$, and then define
$ex_f(n,H)$ as the maximum value of $e_f(G)$ taken over all graphs
with $n$ vertices that do not contain $H$ as a subgraph. Observe
that $ex_f(n,H)=ex(n,H)$ if $f(x)=x/2$, $ex_f(n,H)=ex_p(n,H)$ if
$f(x)=x^p$. Bollob\'as and Nikiforov mentioned that it is important
to study concrete functions. They gave an example
$f(x)=\phi(k)={x\choose k}$, since $\sum_{i=1}^n{d_i\choose k}$
counts the $(k+1)$-vertex subgraphs of $G$ with a dominating vertex.

Denote by $T_r(n)$ the $r$-partite Tur\'an graph of order $n$. In
this paper, using the Bollob\'as--Nikiforov's methods, we give some
results on $ex_{\phi}(n,K_{r+1})$ $(r\geq 2)$ as follows: for
$k=1,2$, $ex_\phi(n,K_{r+1})=e_\phi(T_r(n))$; for each $k$, there
exists a constant $c=c(k)$ such that for every $r\geq c(k)$ and
sufficiently large $n$, $ex_\phi(n,K_{r+1})=e_\phi(T_r(n))$; for a
fixed $(r+1)$-chromatic graph $H$ and every $k$, when $n$ is
sufficiently large, we have
$ex_\phi(n,H)=e_\phi(n,K_{r+1})+o(n^{k+1})$.
\\
[2mm] Keywords: Tur\'an-type problem; degree
sequence; $H$-free\\
[2mm] AMS Subject Classification (2010): 05C35, 05C07.
\end{abstract}

\section{Introduction}
All graphs considered here are finite, undirected, and have no loops
or multiple edges. For the standard graph-theoretic notations the
reader is referred to \cite{B}. Given a graph $G$, let $e(G)$ be the
number of edges of $G$. Denote by $ex(n,H)$ the classical Tur\'an
number, i.e., the maximum number of edges among all graphs with $n$
vertices that do not contain $H$ as a subgraph. Denote by $T_r(n)$
is the $r$-partite Tur\'an graph of order $n$, i.e.,
$ex(n,K_{r+1})=e(T_r(n))$. Given $p\geq 0$ and a graph $G$ whose
degree sequence is $d_1,d_2,\ldots,d_n$, the sum
$e_p(G)=\sum_{i=1}^n d_i^p$ is a much studied parameter in graph
theory, especially for $p=2$. Clearly, $e_1(G)=2e(G)$. In \cite{CY},
Caro and Yuster introduced a Tur\'an-type problem for $e_p(G)$:
given $p\geq 1$, how large can $e_p(G)$ be if $G$ has no subgraph of
a particular type. Denote by $ex_p(n,H)$ be the maximum value of
$e_p(G)$ taken over all graphs with $n$ vertices that do not contain
$H$ as a subgraph. Clearly, $ex_1(n,H)=2ex(n,H)$. It is indeed
interesting to investigate the values of $ex_p(n,H)$ and the
corresponding extremal graphs. It turns out such problems are
usually more difficult for $p>1$ than for $p=1$. Subsequently study
on this topic appeared, such as \cite{BN,BN1,FK,N,P,PT}.

In \cite{CY1}, the authors considered $K_{r+1}$-free graphs and
proved that
\begin{equation}\label{eq1}
ex_p(n,K_{r+1})=e_p(T_r(n))
\end{equation}
for $1\leq p\leq 3$. Therefore, it is interesting to find the values
of $p$ for which equality \eqref{eq1} holds and determine the
asymptotic value of $ex_p(n,K_{r+1})$ for large $n$. In \cite{BN},
Bollob\'as and Nikiforov showed that for every real $p$ ($1\leq p<
r$) and sufficiently large $n$, if $G$ is a graph of order $n$ and
has no clique of order $r+1$, then $ex_p(n,K_{r+1})=e_p(T_r(n))$,
and for every $p\geq r+\lceil\sqrt{2r}\rceil$ and sufficiently large
$n$, $ex_p(n,K_{r+1})>(1+\epsilon)e_p(T_r(n))$ for some positive
$\epsilon=\epsilon(r)$. In \cite{BN1}, Bollob\'as and Nikiforov
proved that if $e_p(G)>(1-1/r)^pn^{p+1}+C$, then $G$ contains more
than $\frac{Cn^{r-p}}{p2^{6r(r+1)+1}r^r}$ cliques of order $r+1$.
Using this statement, they strengthen the Erd\"os--Stone theorem by
using $e_p(G)$ instead of the number of edges.

In \cite{P,PT}, Pikhurko and Taraz discuss even more general
problems. Namely, let $f$ be a non-negative increasing real function
and define $$e_f(G)=\sum_{i=1}^n f(d_i).$$ They generalize the
Tur\'an-type problems and denote by $ex_f(n,H)$ the maximum value of
$e_f(G)$ taken over all graphs with $n$ vertices that do not contain
$H$ as a subgraph. Observe that $ex_f(n,H)=ex(n,H)$ if $f(x)=x/2$,
$ex_f(n,H)=ex_p(n,H)$ if $f(x)=x^p$. They also give several results
for general function $f$. In \cite{BN1}, Bollob\'as and Nikiforov
mentioned that it is important to study concrete functions. They
gave an example $f(x)=\phi(k)={x\choose k}$ and proposed an open
problem, since the combinatorial implications of this question are
obvious, i.e., $\sum_{i=1}^n{d_i\choose k}$ counts the
$(k+1)$-vertex subgraphs of $G$ with a dominating vertex.

In this paper, we will study this new Tur\'an-type problem on degree
sequence. Using the Bollob\'as--Nikiforov's methods and according
their results \cite{BN}, we give some results on
$ex_{\phi}(n,K_{r+1})$, i.e., the maximum value of
$$e_\phi(G)=\sum_{i=1}^n{d_i\choose k}$$ taken over all
$K_{r+1}$-free graphs with $n$ vertices. We assume $k\leq n/2$ in
this paper.

Our first result shows that for $k=1,2$, the extremal graph that
yields $ex_\phi(n,K_{r+1})$ is exactly the Tur\'an graph $T_r(n)$.

\begin{thm}\label{thm1}
Let $r\geq 2$ be a positive integer and $k=1,2$. Then
\begin{equation}\label{eq2}
ex_\phi(n,K_{r+1})=e_\phi(T_r(n)),
\end{equation}
where $T_r(n)$ is the $r$-partite Tur\'an graph of order $n$.
\end{thm}

Actually, Theorem \ref{thm1} is not true for the case of $k\geq 3$.
We consider the triangle-free graphs. The structure of triangle-free
graphs such that $e_\phi(G)=ex_\phi(n,K_{3})$ for $k=3$ are
characterized. Using the similar method, we also characterize the
structure of triangle-free graphs such that $e_p(G)=ex_p(n,K_{3})$
for $p\geq 4$.

Similarly, a natural problem arises: given $r\geq 2$, determine the
integers $k$, for which equality \eqref{eq2} holds. We will deduce
the following result.
\begin{thm}\label{thm2}
For each $k$, there exists a constant $c=c(k)$ such that for every
integer $r\geq c(k)$ and sufficiently large $n$,
$$ex_\phi(n,K_{r+1})=e_\phi(T_r(n)),$$
where $c(k)=1+\lim\limits_{n\rightarrow \infty}(n/c_1)$ and $0\leq
c_1\leq n-k$ is the root of the function
$1-x\Psi(n-x+1)+x\Psi(n-x-k+1)$, where $\Psi(x)$ is the Digamma
function.
\end{thm}

Computer calculations show that $c(k)=5/2$ when $k=2$, and $c(k)=5$
when $k=3$.

For a fixed $(r+1)$-chromatic graph $H$, we show that for every
$r\geq 2$ and $k$, we have
$ex_\phi(n,H)=e_\phi(n,K_{r+1})+o(n^{k+1})$.
\begin{thm}\label{thm3}
Let $H$ be a graph of order $n$ with $\chi(H)=r+1$ $(r\geq 2)$. When
$n$ is sufficiently large, we have
$$ex_\phi(n,H)=e_\phi(n,K_{r+1})+o(n^{k+1}).$$
\end{thm}

\section{Proof of Theorem \ref{thm1}}
In order to prove Theorem \ref{thm1}, we need the following theorem
of Erd\"os \cite{E}, which characterizes the maximal degree
sequences of graphs without a $K_{r+1}$.

\begin{lem}[\cite{E}]\label{lem1}
Let $G=(V,E)$ be a graph without a $K_{r+1}$. Then there is a
$r$-partite graph $G'=(V,E')$ such that for every $v\in V$,
$d_G(v)\leq d_{G'}(v)$.
\end{lem}

If $G$ and $G'$ are as in Lemma \ref{lem1}, then for every integer
$k$, $$\sum_{i=1}^k{d_G(v)\choose k}\leq
\sum_{i=1}^k{d_{G'}(v)\choose k},$$ i.e., $e_\phi(G)\leq
e_\phi(G')$. Thus, the following corollary immediately follows.
\begin{cor}\label{cor1}
For every $n\geq r\geq 1$ and integer $k$, there exists a complete
$r$-partite graph $G$ with $n$ vertices such that
$ex_\phi(n,K_{r+1})=e_\phi(G)$.
\end{cor}

Actually, Pikhurko \cite{P} obtained the following general result.
\begin{cor}
For every $n\geq r$ and non-decreasing function $f:
\mathbb{N}\rightarrow \mathbb{R}$, there exists a complete
$r$-partite graph $G$ with $n$ vertices such that
$ex_f(n,K_{r+1})=e_f(G)$.
\end{cor}

\noindent{\bf Proof of Theorem \ref{thm1}.} We only need to consider
the case of $k=2$, since $ex_\phi(n,H)=2ex(n,H)$ for $k=1$. Let $G'$
be any complete $r$-partite graph that has at least two vertex
classes $X$ and $Y$ with $|X|=x$, $|Y|=y$ and $x-y>1$. The complete
$r$-partite graph $G''$ is obtained from $G'$ by transferring a
vertex from $X$ to $Y$. Then we only need to prove
$e_\phi(G'')>e_\phi(G')$ for $k=2$. Observe that $n\geq x+y$,
$x-y-1>0$, and
\begin{eqnarray*}
&&e_\phi(G'')-e_\phi(G')\\[2mm]
&=&(y+1){n-y-1\choose 2}+(x-1){n-x+1\choose 2} -y{n-y\choose
2}-x{n-x\choose 2}.
\end{eqnarray*}
Let $g(n,x,y)=e_\phi(G'')-e_\phi(G')$. By some elementary
calculations, we have that $g(n,x,y)$ is monotonously increasing on
$n$.
Therefore, $g(n,x,y)\geq g(x+y,x,y)>g(2y+1,y+1,y)>0$, since
$g(x+y,x,y)$ is monotonously increasing on $x$.\qed\\

Observe that Theorem \ref{thm1} is not true for $k\geq 3$. This can
be seen from the fact that for some $3\leq k\leq n/2$,
$e_\phi(K_{\lfloor n/2-1\rfloor,\lceil
n/2+1\rceil})>e_\phi(T_3(n))$. In the following, we will consider
triangle-free graphs for $k=3$. Let $G^*$ be the triangle-free graph
such that $e_\phi(G^*)=ex_\phi(n,K_{3})$. From Corollary \ref{cor1},
$G^*$ is a complete bipartite graph and we assume $G=K_{x,n-x}$. We
start by computing $e_\phi(G)$ for $k=3$. Set $e=e(G)=x(n-x)$.
\begin{eqnarray*}
e_\phi(G)&=&\frac{x(n-x)(n-x-1)(n-x-2)+x(n-x)(x-1)(x-2)}{6}\\
&=&\frac e 6(-2e+n^2-3n+4)\\
&=&\frac 1 3\left(\frac{n^2-3n+4}{4}\right)^2-\frac 1
3\left(\frac{n^2-3n+4}{4}-e\right)^2.
\end{eqnarray*}
Let $\mathcal{G}_n^*$ denote the set of (unbalanced) complete
bipartite graphs $K_{x,n-x}$ for which $e=x(n-x)$ is closest to
$(n^2-3n+4)/4$. Observe that $\mathcal{G}_n^*$ usually consists of a
single graph, but an elementary number of theoretical argument shows
that there are also infinitely many cases when
$\mathcal{G}_n^*=\{K_{x,\ n-x},K_{x+1,\ n-x-1}\}$. Indeed, in these
cases, $n,\, k$ satisfying that
$$\frac{n^2-3n+4}{4}-x(n-x)=(x+1)(n-x-1)-\frac{n^2-3n+4}{4},$$
i.e., $4xn-4x^2-4x-n^2+5n=6$. Computer calculations show that the
first few cases are
\begin{eqnarray*}
(n,k)=&&(7,1),\ (10, 2),\ (18,5),\ (23,7),\ (35,12),\ (42,15),\
(58,22),\ (67,26),\\
&&(87, 35),\ (98,40),\ (122,51),\ (135,57),\ (163,70),\ (178,77),\
\ldots.
\end{eqnarray*}
We conclude this result as follows.
\begin{thm}
If $G$ is a triangle-free graph of order $n$ and
$e_\phi(G)=ex_\phi(n,K_{3})$ for $k=3$, then $G\in \mathcal{G}_n^*$,
where $\mathcal{G}_n^*$ denotes the set of (unbalanced) complete
bipartite graphs $K_{x,n-x}$ for which $e=x(n-x)$ is closest to
$(n^2-3n+4)/4$.\qed
\end{thm}

Recall that in \cite{BN}, Bollob\'as and Nikiforov proved that for
$0< p\leq 3$, $ex_p(n,K_3)=e_p(T_2(n))$, and showed that for every
$\varepsilon>0$, there exists a $\delta$ such that if $p>3+\delta$
then $ex_p(n,K_3)>(1+\varepsilon)e_p(T_2(n))$ for sufficiently large
$n$. They also proved that if $ex_p(n,K_{r+1})=e_p(T_r(n))$ for some
$K_{r+1}$-free graph $G$ of order $n$, then $G$ is a complete
$r$-partite graph having $r-1$ vertex classes of size $cn+o(n)$,
where $c=c(p,r)$ is a constant. Using the above method, we will show
that among all triangle-free graphs, some unbalanced complete
bipartite graphs will attain the maximum value of $e_p(G)$ for
$p\geq 4$. Let $G^*$ be the triangle-free graph such that
$e_p(G^*)=ex_p(n,K_{3})$. Similarly, $G^*$ must be a complete
bipartite graph and we assume $e(G^*)=a n^2$, where $a=a(p)$ is a
constant. Let $G=K_{x,n-x}$. Set $e=e(G^*)=x(n-x)$ and then we can
assume $x=\frac {n-\sqrt{n^2-4e}}{2}$, $n-x=\frac
{n+\sqrt{n^2-4e}}{2}$. Let
\begin{eqnarray*}
h(n,e)&=&e_p(G)=x(n-x)^p+(n-x)x^p=e\left((n-x)^{p-1}+x^{p-1}\right)\\[2mm]
&=&e\left(\left(\frac {n+\sqrt{n^2-4e}}{2}\right)^{p-1}+\left(\frac
{n-\sqrt{n^2-4e}}{2}\right)^{p-1}\right).
\end{eqnarray*}

For $p=4$, we have $h(n,e)=\frac {n^5}{12}-3n\left(e-\frac {n^2}
6\right)^2$, which implies that $a=1/6\approx 0.166666666$ for
$p=4$. For $p=5$, we have $h(n,e)=2e^3-4n^2e^2+n^4e$, which implies
that $a=\frac 2 3-\frac{\sqrt{10}}{6}\approx 0.139620390$ for $p=5$.
Computer calculations show that the value of $a=\frac 1
3-\frac{\sqrt{10}}{15}\approx 0.122514822$ for $p=6$, $a\approx
0.1093797828$ for $p=7$, $a\approx 0.09876572056$ for $p=8$,...
Observe that the value of $a$ will decrease when $p$ increases.

Actually, the maximum value of $e_p(G)$ will be attained when $e=a
n^2$ is the root of $\frac{\partial(h(n,e))}{\partial e}$, i.e.,
$e=a n^2$ is the root of the following equation:
$$\frac{n^2-2e(p+1)+n\sqrt{n^2-4e}}{n^2-2e(p+1)-n\sqrt{n^2-4e}}
=\left(\frac{n-\sqrt{n^2-4e}}{n+\sqrt{n^2-4e}}\right)^{p-2}.$$
Obviously, the right side of the above equality is positive, and so
is the left side. Then we have $n^2-2e(p+1)-n\sqrt{n^2-4e}>0$ or
$n^2-2e(p+1)+n\sqrt{n^2-4e}<0$. It implies that
$$\left(n^2-2e(p+1)\right)^2>\left(n\sqrt{n^2-4e}\right)^2,$$
from which we obtain that $e>\frac p {(p+1)^2}n^2$. It is very
surprise that the value of $\frac p {(p+1)^2}n^2$ is very close to
what we need. Some calculations show that the value of $\frac p
{(p+1)^2}$ is $0.16$ for $p=4$, $0.1388888889$ for $p=5$,
$0.1224489796$ for $p=6$, $0.1093750000$ for $p=7$, $0.09876543210$
for $p=8$,..., which is very close to the exact value of $a$,
especially for $p\geq 7$.

Therefore, we can conclude as follows.
\begin{thm} If $G$ is a triangle-free graph of order
$n$ and $e_p(G)=ex_p(n,K_{3})$ for $p\geq 7$, then $G\in
\mathcal{G}_n^{**}$, where $\mathcal{G}_n^{**}$ denotes the set of
(unbalanced) complete bipartite graphs $K_{x,n-x}$ for which
$e=x(n-x)$ is closest to $\frac p {(p+1)^2}n^2$.\qed
\end{thm}

\section{Proof of Theorem \ref{thm2}}

\noindent{\bf Proof of Theorem \ref{thm2}.} From Corollary
\ref{cor1}, we know that for $K_{r+1}$-free graphs $G$ of order $n$,
if $e_\phi(G)$ attains a maximum, then $G$ is a complete $k$-partite
graph. Let $G$ be any complete $k$-partite graph. Notice that every
complete $k$-partite graph can be determined uniquely by the size of
its vertex classes. So we assume the sizes of vertex classes of $G$
are $n_1,n_2,\ldots,n_r$, satisfying that $1\leq n_1\leq n_2\leq
\cdots\leq n_r$. Therefore, $e_\phi(n,K_{r+1})$ equals to
$$\max\left\{\sum_{i=1}^rn_i{n-n_i\choose k}: n_1+n_2+
\cdots+ n_r=n,\ 1\leq n_1\leq n_2\leq \cdots\leq n_r\right\}.$$

Let $g(x,n,k)=x{n-x\choose k}$ and we have
$$\frac {\partial
g(x,n,k)}{\partial x}=x{n-x\choose
k}\left(1-x\Psi(n-x+1)+x\Psi(n-x-k+1)\right),$$ where $\Psi(x)$ is
the Digamma function.  Routine calculations show that the function
$x{n-x\choose k}$ increases for $0\leq x\leq c_1$, decreases for
$c_1\leq x\leq n-k$, and is concave for $0\leq x\leq c_2$, where
$c_1$ and $c_2$ denote the root of $\frac {\partial
g(x,n,k)}{\partial x}$ and $\frac {\partial^2 g(x,n,k)}{\partial
x^2}$, respectively.

Now suppose $G^*$ is the complete $k$-partite graph which attains
the value of $e_\phi(n,K_{r+1})$. Denote by
$n^*_1,n^*_2,\ldots,n^*_r$ the sizes of vertex classes of $G^*$.

If $n^*_r\leq c_2$, then the concavity of $x{n-x\choose k}$ implies
that $n^*_r-n^*_1\leq 1$, and the proof is thus completed.

Now we will assume $n^*_r>c_2$. In this case, we claim that
$n^*_1\geq c_1$. Since otherwise, adding $1$ to $n^*_r$ and
subtracting $1$ from $n^*_1$, the value $x{n-x\choose k}$ will
increase, contradicting to the choice of $G^*$. Hence, we have
$$c_1\leq n^*_1\leq \frac n r \leq \frac n {c(k)}=\frac{n}{1+\lim\limits_{n\rightarrow \infty}(n/c_1)},$$
i.e., $1+\lim\limits_{n\rightarrow \infty}(n/c_1)\leq
\frac{n}{c_1}$, a contradiction when $n\rightarrow \infty$. \qed

\section{Proof of Theorem \ref{thm3}}
To prove Theorem \ref{thm3}, we need the following result (for a
proof see, e.g., \cite{B}, Theorem 33, p. 132).

\begin{lem}[\cite{B}]\label{lem2}
Suppose $H$ is an $(r+1)$-chromatic graph. Every $H$-free graph $G$
of sufficiently large order $n$ can be made $K_{r+1}$-free by
removing $o(n^2)$ edges.
\end{lem}

\noindent{\bf Proof of Theorem \ref{thm3}.} Let $G$ be a
$K_{r+1}$-free graph of order $n$ such that
$ex_\phi(n,K_{r+1})=ex_\phi(G)$. By Corollary \ref{cor1}, $G$ is
$r$-partite, which implies that $G$ is $H$-free. Thus, we have
$ex_\phi(n,H)\geq ex_\phi(G)=e_\phi(n,K_{r+1})$. Now let $G'=(V,E')$
be an $H$-free graph of order $n$ with $ex_\phi(n,H)=ex_\phi(G')$.
From Lemma \ref{lem2}, there exists a $K_{r+1}$-free graph
$G''=(V,E'')$ that may be obtained from $G'$ by removing at most
$o(n^2)$ edges.

Denote by $d'_i$ and $d''_i$ the degree of vertex $v_i$ in $G'$ and
$G''$, respectively. Obviously, $d'_i\geq d''_i$. For each $1\leq
i\leq n$, we consider the difference ${d'_i\choose k}-{d''_i\choose
k}$. Using the elementary calculus, we have $${d'_i\choose
k}-{d''_i\choose k}<k(d'_i-d''_i)\left(\frac e
k\right)^k(d'_i)^{k-1}<k(d'_i-d''_i)\left(\frac e
k\right)^kn^{k-1}.$$ Summing this inequality for all $1\leq i\leq
n$, we have
$$ex_\phi(G')-ex_\phi(G'')<k\left(\frac e
k\right)^kn^{k-1}\cdot 2\left(e(G')-e(G'')\right)=o(n^{k+1}).$$ The
proof is thus completed. \qed

\end{document}